\documentclass[a4paper,UKenglish]{lipics_md}
\usepackage{microtype}
\usepackage{graphicx}
\usepackage{tikz}

\title{Tri-Partitions and Bases of an Ordered Complex\footnote{
        This project has received funding from the European Research Council (ERC)
        under the European Union's Horizon 2020 research and innovation programme
        (grant agreement No 78818 Alpha).
	It is also partially supported by the DFG Collaborative Research Center TRR 109,
	`Discretization in Geometry and Dynamics',
	through grant no.\ I02979-N35 of the Austrian Science Fund (FWF).}}

\author[1]{Herbert Edelsbrunner}
\author[1]{Katharina \"{O}lsb\"{o}ck}
\affil[1]{IST Austria (Institute of Science and Technology Austria),
	Klosterneuburg, \\
        Austria, \texttt{edels@ist.ac.at}, \texttt{katharina.oelsboeck@ist.ac.at}}
\authorrunning{H. Edelsbrunner and K. \"{O}lsb\"{o}ck}


\newcommand {\mm}[1] {\ifmmode{#1}\else{\mbox{\(#1\)}}\fi}

\newcommand{\denselist}{\itemsep 0pt\parsep=1pt\partopsep 0pt}

\newcommand{\ignore}[1]{}

\newcommand{\ourproof}{\begin{proof}}
\newcommand{\eop}{\end{proof}}

\newcommand{\Zspace}       {\mm{{\mathbb Z}}}
\newcommand{\Acal}         {\mm{{\mathcal A}}}
\newcommand{\Ecal}         {\mm{{\mathcal E}}}
\newcommand{\Fcal}         {\mm{{\mathcal F}}}

\newcommand{\Cycle}[2]     {\mm{{\sf z}_{#1}{({#2})}}}
\newcommand{\Cocycle}[2]   {\mm{{\sf z}^{#1}{({#2})}}}
\newcommand{\Chain}[2]     {\mm{{\sf c}_{#1}{({#2})}}}
\newcommand{\Xain}[2]      {\mm{{\sf x}_{#1}{({#2})}}}
\newcommand{\Cochain}[2]   {\mm{{\sf c}^{#1}{({#2})}}}
\newcommand{\Coxain}[2]    {\mm{{\sf x}^{#1}{({#2})}}}

\newcommand{\Bgr}[1]       {\mm{{\sf B}_{#1}}}
\newcommand{\coBgr}[1]     {\mm{{\sf B}^{#1}}}
\newcommand{\Cgr}[1]       {\mm{{\sf C}_{#1}}}
\newcommand{\coCgr}[1]     {\mm{{\sf C}^{#1}}}
\newcommand{\rHgr}[1]      {\mm{\tilde{\sf H}_{#1}}}
\newcommand{\corHgr}[1]    {\mm{\tilde{\sf H}^{#1}}}
\newcommand{\Zgr}[1]       {\mm{{\sf Z}_{#1}}}
\newcommand{\coZgr}[1]     {\mm{{\sf Z}^{#1}}}
\newcommand{\bgr}[1]       {\mm{{\sf b}_{#1}}}
\newcommand{\cgr}[1]       {\mm{{\sf c}_{#1}}}
\newcommand{\xgr}[1]       {\mm{{\sf x}_{#1}}}
\newcommand{\zgr}[1]       {\mm{{\sf z}_{#1}}}
\newcommand{\cocgr}[1]     {\mm{{\sf c}^{#1}}}
\newcommand{\coxgr}[1]     {\mm{{\sf x}^{#1}}}
\newcommand{\cozgr}[1]     {\mm{{\sf z}^{#1}}}

\newcommand{\Skeleton}[2]  {\mm{{#1}^{({#2})}}}

\newcommand{\rBetti}[1]    {\mm{{\tilde{\beta}}_{#1}}}
\newcommand{\corBetti}[1]  {\mm{{\tilde{\beta}}^{#1}}}

\newcommand{\rEuler}       {\mm{{\tilde{\chi}}}}
\newcommand{\nSx}[1]       {\mm{n_{#1}}}
\newcommand{\nSxB}[1]      {\mm{n_{#1}^\circ}}
\newcommand{\nSxD}[1]      {\mm{n_{#1}^\bullet}}
\newcommand{\ncoSx}[1]     {\mm{n_{#1}}}
\newcommand{\ncoSxB}[1]    {\mm{n^{#1}_\circ}}
\newcommand{\ncoSxD}[1]    {\mm{n^{#1}_\bullet}}

\newcommand{\Low}[1]       {\mm{\rm low}{({#1})}}
\newcommand{\Left}[1]      {\mm{\rm left}{({#1})}}
\newcommand{\rank}[1]      {\mm{\rm rank\,}{#1}}
\newcommand{\dime}[1]      {\mm{\rm dim\,}{#1}}

\newcommand{\ssx}          {\mm{\sigma}}

\newcommand{\card}[1]      {|{#1}|}

\newcommand{\ourparagraph}[1] {\vspace{0.1in} \noindent \textbf{#1}}
\newcommand{\Skip}[1]      {}
\definecolor{blue-green}{rgb}{0.0, 0.87, 0.87}

\begin{document}
\maketitle

\begin{abstract}
  Generalizing the decomposition of a connected planar graph into a tree and a
  dual tree, we prove a combinatorial analog of the classic
  Helmholtz--Hodge decomposition of a smooth vector field.
  Specifically, we show that for every polyhedral complex, $K$, and every dimension, $p$,
  there is a partition of the set of $p$-cells into a maximal $p$-tree,
  a maximal $p$-cotree, and a collection of $p$-cells whose cardinality
  is the $p$-th reduced Betti number of $K$.
  Given an ordering of the $p$-cells, this tri-partition is unique,
  and it can be computed by a matrix reduction algorithm that
  also constructs canonical bases of cycle and boundary groups.
\end{abstract}

\section{Introduction}
\label{sec:1}

Given a connected graph embedded on the sphere, it is well known that we can
split the graph into a spanning tree and a dual tree whose nodes
are the faces.
This is best visualized by rotating each edge that is not in the spanning tree,
making sure they meet at points chosen to represent the faces.
This split is similar in spirit to the Helmholtz decomposition of a smooth vector field
on the sphere into a rotation-free component and a divergence-free component \cite{Hel58}.
If the graph is embedded on a surface with non-zero genus, then
the split does not exhaust all edges and
the unused ones correspond to the third, harmonic component
of the Helmholtz--Hodge decomposition on this surface \cite{Hod41}.

Thinking of this split as a theorem about the edges of a connected planar graph,
we are interested in its generalization to complexes and to cells
of any dimension.
Such a generalization promises a geometric interpretation
of algebraic concepts in homology and their relations.
Beyond this theoretical interest in the structure of complexes,
we are motivated by geometric modeling tasks in which holes are of
central importance.
An example are cell membrane proteins with functional channels for ion transport.
We believe that our structural results can be helpful in discovering
and manipulating hole systems,
but this is the topic of future work.

\begin{figure}[hbt]
  \centering
  \resizebox{!}{2.4in}{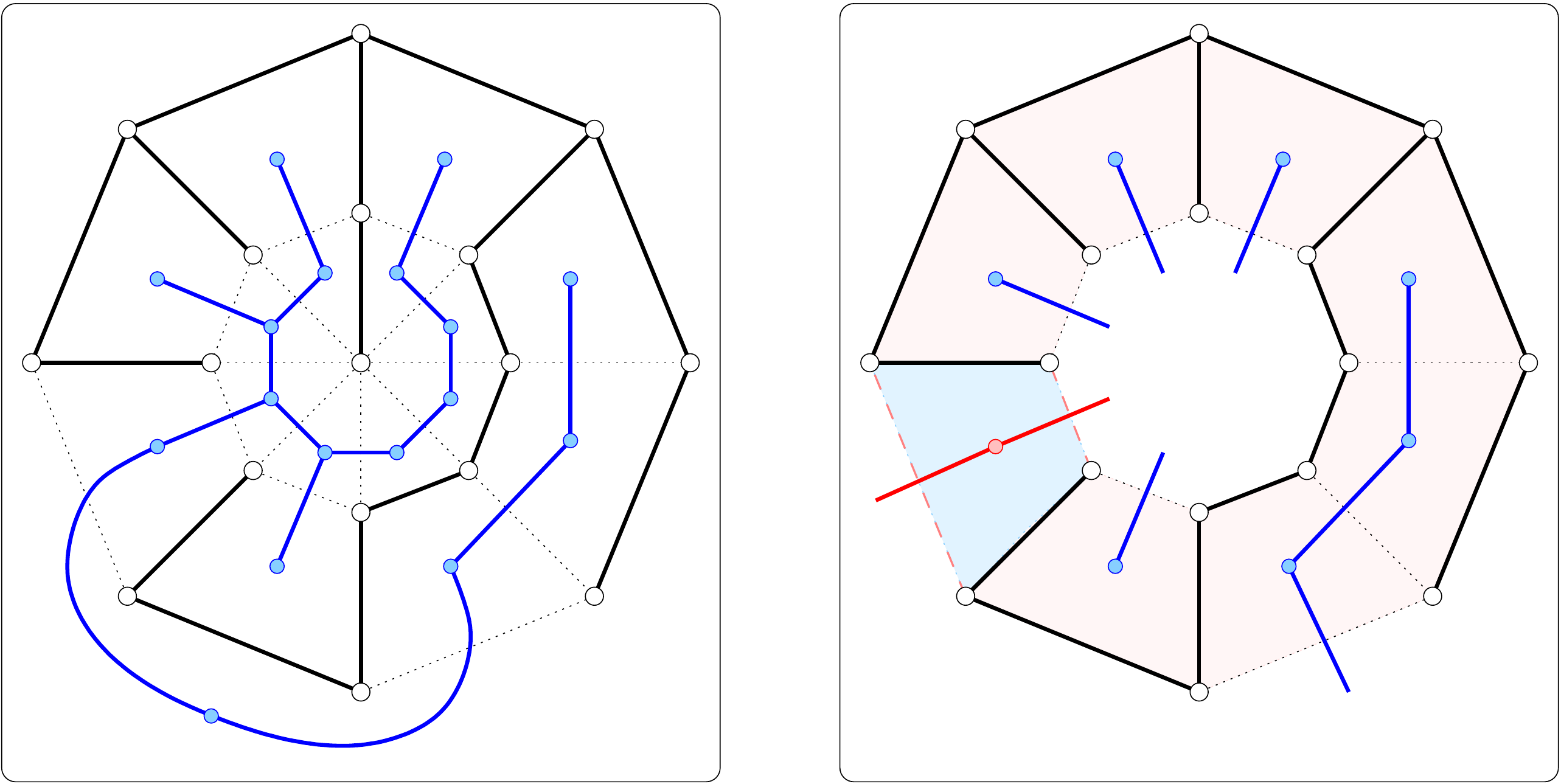}
  \caption{\emph{Left panel:} a graph of solid and dotted black edges
    embedded in the plane.
    The solid edges form a spanning tree of the graph, and the dotted
    edges intersect the dual blue edges
    that form a spanning tree of the dual graph.
    \emph{Right panel:} an annulus decomposed into eight quadrangles.
    Besides a solid black spanning tree and a dotted cotree,
    we show the two dashed edges of a cocycle.
    To get a maximal cotree, we add one of the dashed edges to the cotree.}
  \label{fig:disk}
\end{figure}

\ourparagraph{Two examples.}
Our results are combinatorial and algorithmic.
To get a first impression, consider the planar graph forming a wheel of
$n_0 = 17$ vertices and $n_1 = 32$ edges drawn in the left panel of
Figure \ref{fig:disk}.
Every spanning tree consists of $n_0 - 1 = 16$ edges,
and if we interpret this tree as a barrier between the $n_2 = 17$
$2$-dimensional regions defined by the embedding of the graph,
then each non-tree edge splits the regions into two connected collections.
In other words, the regions form a dual tree, drawn with blue edges
crossing the dotted non-tree edges.

As already noted in \cite{Big71,RoRe78}, the bi-partition of the edges
is best generalized to a tri-partition if the graph is embedded on a closed
surface with positive genus.
To further free ourselves from the implicit definition of regions,
we consider complexes in which the $2$-dimensional cells are explicitly
specified, so it no longer matters where the complex is embedded.
An example is the complex of $n_0 = 16$ vertices, $n_1 = 24$ edges,
and $n_2 = 8$ quadrangles drawn in the right panel of Figure \ref{fig:disk}.
An unpleasant consequence is that the dual is no longer necessarily
closed as it may have edges with missing endpoints.
The appropriate formalism is therefore cohomology,
where we worry about edges and their incident regions rather than their
incident vertices.
In this formalism, a \emph{cocycle} is a set of edges such that every region
in the complex is incident to an even number of edges in this set,
and a \emph{cotree} is a set of edges that does not contain a cocycle.
Returning to the annulus in Figure \ref{fig:disk}, we notice that the spanning tree
cannot go completely around the `hole', so we get a cocycle that connects
the outer boundary with the inner boundary.
In other words,
the complement of the spanning tree is a cotree-with-a-single-cocycle.
Alternatively, we could add one edge of the cocycle to the tree to get
a tree-with-a-single-cycle whose complement is a genuine cotree.
Better yet, we partition the edges into a maximal tree, a maximal cotree,
and one leftover edge.
In the general case, there can be more leftover edges, namely one for each `hole'.

\ourparagraph{Results and prior work.}
The partition of the edge set has been studied for connected
graphs embedded on orientable closed surfaces.
Biggs proved in 1971 that such a graph splits into a spanning tree
and a complementary subgraph of the dual that contains a spanning tree
of the dual \cite{Big71}.
Rosenstiehl and Read sharpened the result in 1978 by observing
that the complementary subgraph of the dual splits into a
(dual) spanning tree and $2g$ additional edges,
in which $g$ is the genus of the surface \cite{RoRe78}.
The first result of this paper generalizes this split to
complexes and to cells beyond edges.
Specifically, we prove that for every polyhedral complex, $K$, and for every
dimension, $p$, the set of $p$-cells can be partitioned into a maximal $p$-tree,
a maximal $p$-cotree, and a set of leftover $p$-cells whose
cardinality is the $p$-th reduced Betti number of $K$.
These properties suggest that the tri-partition is a discrete analog of the
Helmholtz--Hodge decomposition of a smooth vector field into
a rotation-free component, a divergence-free component, and
a harmonic component whose dimension is the rank of the homology group
\cite{Hel58,Hod41} but see also \cite{BNPB13,PoPr00}.
An algebraic analog of this decomposition was introduced by Eckmann in 1945,
namely that the $p$-th chain group satisfies
$\Cgr{p} \simeq \Bgr{p} \oplus \coBgr{p} \oplus \rHgr{p}$ \cite{Eck45}
but see also \cite{Dod74,Fri98}.

The tools for establishing the tri-partition are the matrix reduction
algorithms developed in the context of persistent homology \cite{EdHa10},
the theorem on the invariance of birth-death pairs proved in \cite{CEM06},
and the duality between homology and relative cohomology noted in \cite{dSMV11}.
Importantly, the tri-partition implies canonical bases of the cycle,
boundary, and homology groups as well as of their counterparts in cohomology.
More specifically, for each monotonic ordering of a polyhedral complex,
there is a unique tri-partition and a unique collection of bases.

\ourparagraph{Outline.}
Section \ref{sec:2} provides background from algebraic topology,
including algorithms for the ranks of homology and cohomology groups.
Section \ref{sec:3} introduces the tri-partition of a polyhedral complex.
Section \ref{sec:4} describes the related bases
in homology and cohomology and proves their properties.
Section \ref{sec:5} concludes this paper.
Appendix \ref{app:A} proves that the tri-partitions form matroids.

\section{Background}
\label{sec:2}

We will make frequent use of homology and cohomology groups;
see \cite{Hat02,Mun84} for general background on these topics.
To keep the discussion elementary, we use $\Zspace / 2 \Zspace$ coefficients
so that cycles and cocycles can be treated as sets.

\ourparagraph{Polyhedral complexes.}
A \emph{$p$-cell}, $\ssx$, is a $p$-dimensional convex polytope,
and we write $\dime{\ssx} = p$ for its \emph{dimension}.
A hyperplane \emph{supports} $\ssx$ if it has a non-empty intersection
with the polytope and the polytope is contained in a closed half-space
bounded by the hyperplane.
A \emph{face} of $\ssx$ is the intersection with a supporting hyperplane;
it is a convex polytope of dimension at most $p$.
We call $\ssx$ a \emph{coface} of its faces.
A \emph{polyhedral complex}, $K$, is a collection of cells that is closed
under taking faces such that the intersection of any two cells is a face of both.
By convention, we require that the empty cell is part of $K$;
its dimension is~$-1$ and it is a face of every cell.
A cell is \emph{maximal} if it has no proper coface in $K$.
The \emph{dimension} of $K$ is the maximum dimension of any of its cells.
The \emph{$p$-skeleton} contains all cells of dimension $p$ or less
and is denoted $\Skeleton{K}{p}$.
We write $K^p = \Skeleton{K}{p} \setminus \Skeleton{K}{p-1}$
for the set of $p$-cells in $K$, and $\nSx{p} = \card{K^p}$ for its cardinality,
noting that $\nSx{p} = 0$ for $p$ smaller than $-1$ and larger than $\dime{K}$.
The \emph{Euler characteristic} is the alternating sum of cell numbers,
and since the empty cell is included, we decorate it with a tilde:
$\rEuler = \sum_p (-1)^p \nSx{p}$.
The Euler--Poincar\'{e} Formula asserts that the Euler characteristic
is the alternating sum of Betti numbers and therefore a topological invariant.
We formally state the result now and provide the definition of
the Betti numbers later.
\begin{proposition}[Euler--Poincar\'{e}]
  \label{prop:EulerPoincare}
  Every polyhedral complex satisfies
  $\rEuler = \sum_p (-1)^p \rBetti{p}$.
\end{proposition}
To represent a polyhedral complex in the computer, it is common to
order the cells --- arbitrarily or otherwise --- and to store the
face relation in matrix form.
Letting $\ssx_0, \ssx_1, \ldots, \ssx_m$ be the ordering,
the \emph{boundary matrix}, $\partial [0..m, 0..m]$, is defined by
\begin{align}
  \!\partial [i,j]  &=  \left\{ \begin{array}{cl}
                        \!1  &  \mbox{\rm if~} \ssx_i \subseteq \ssx_j
                              \mbox{\rm ~and~} \dime{\ssx_i} = \dime{\ssx_j}-1 , \\
                        \!0  &  \mbox{\rm otherwise}.
                      \end{array} \right.
\end{align}
In words: column $j$ of $\partial$ stores the codimension $1$ faces of $\ssx_j$
and row $i$ stores the codimension~$1$ cofaces of $\ssx_i$.
Throughout this paper, we use \emph{monotonic orderings} in which every cell
is preceded by its faces.
The boundary matrix of a monotonically ordered polyhedral complex is upper-triangular.
A \emph{filtration} of $K$ is a nested sequence of subcomplexes that ends with $K$.
An example are the prefixes of a monotonic ordering:
$K_\ell = \{ \ssx_0, \ssx_1, \ldots, \ssx_\ell \}$ is a complex,
for every $0 \leq \ell \leq m$, 
and $K_0 \subseteq K_1 \subseteq \ldots \subseteq K_m$ is a filtration of $K$.

\ourparagraph{Homology.}
Since we use $\Zspace / 2 \Zspace$ coefficients, we define a \emph{$p$-chain}
as a subset of the $p$-cells, $\cgr{p} \subseteq K^p$.
Accordingly, the \emph{sum} of two $p$-chains is their symmetric difference,
and this operation defines a group, denoted $\Cgr{p} (K)$.
The \emph{boundary} of a $p$-chain is the $(p-1)$-chain, $\partial \cgr{p}$,
that consists of all $(p-1)$-cells shared by an odd number
of $p$-cells in $\cgr{p}$.
A \emph{$p$-cycle} is a $p$-chain with empty boundary,
and a \emph{$p$-boundary} is the boundary of a $(p+1)$-chain.
Since we include the empty cell, the boundary of a vertex is
this empty cell and therefore not empty.
In contrast to conventional homology theory,
a single vertex is therefore not a $0$-cycle,
but a pair of vertices is.
The $p$-boundaries and $p$-cycles form subgroups of $\Cgr{p}$,
and because taking the boundary twice always gives the empty set,
the former is a subgroup of the latter:
$\Bgr{p} (K) \subseteq \Zgr{p} (K) \subseteq \Cgr{p} (K)$.
A $p$-cycle in $\Bgr{p} (K)$ is sometimes referred to as \emph{trivial},
and two $p$-cycles are \emph{homologous} if they differ by a $p$-boundary.
The \emph{$p$-th (reduced) homology group} consists of all classes
of homologous $p$-cycles:
$\rHgr{p} (K) = \Zgr{p} (K) / \Bgr{p} (K)$.
The \emph{$p$-th (reduced) Betti number} is the rank of the $p$-th homology group.
Since we use modulo-$2$ arithmetic, this rank
is the binary logarithm of the cardinality, and we write
$\rBetti{p} = \rBetti{p} (K) = \rank{\rHgr{p} (K)} = \log_2 \card{\rHgr{p} (K)}$.
We call $K$ \emph{acyclic} if all Betti numbers vanish.
The smallest polyhedral complex is $K = \{ \emptyset \}$, with $\rBetti{-1} = 1$
and all other Betti numbers zero.
The smallest acyclic polyhedral complex consists of a vertex and the empty cell.

\ourparagraph{Cohomology.}
While we collect information about the complex using homology,
we collect information about the complement using relative cohomology;
see also \cite{dSMV11}.
We begin with the definitions for cohomology.
A \emph{$p$-cochain} is a subset of the $p$-cells, $\cocgr{p} \subseteq K^p$.
Its \emph{coboundary}, $\delta \cocgr{p}$, consists of all $(p+1)$-cells
that have an odd number of faces in $\cocgr{p}$.
A \emph{$p$-cocycle} is a $p$-cochain with empty coboundary,
and a \emph{$p$-coboundary} is the coboundary of a $(p-1)$-cochain.
Again we get groups, $\coBgr{p} (K) \subseteq \coZgr{p} (K) \subseteq \coCgr{p} (K)$,
which we distinguish from the boundary, cycle, and chain groups
by writing the dimension as superscript.
The \emph{$p$-th (reduced) cohomology group} consists of all classes
of cohomologous $p$-cocycles: $\corHgr{p} (K) = \coZgr{p} (K) / \coBgr{p} (K)$.
We write
$\corBetti{p} = \corBetti{p} (K)
              = \rank{\corHgr{p} (K)} = \log_2 \card{\corHgr{p} (K)}$.
For example, if $K = \{\emptyset\}$, then $\corBetti{-1} = 1$ and $\corBetti{p} = 0$
for all $p \geq 0$.
We will see shortly that $\corBetti{p} = \rBetti{p}$ for all $p$,
so $K$ is acyclic iff $\corBetti{p} = 0$ for all dimensions $p$.
As a general intuition, $\corBetti{p}$ is the number of cuts
needed to remove all non-trivial $p$-th cohomology.

Relative cohomology is similar but defined for a pair, $(K, L)$,
in which $L$ is a subcomplex of $K$.
The \emph{relative $p$-cochains} are the $p$-cochains in $K \setminus L$,
and we notice that their coboundaries are also in $K \setminus L$.
We therefore define the \emph{relative $p$-cocycles} as the $p$-cocycles
in $K \setminus L$,
and the \emph{relative $p$-coboundaries} as the $p$-coboundaries in $K \setminus L$.
As before, we get three nested groups,
$\coBgr{p} (K,L) \subseteq \coZgr{p} (K,L) \subseteq \coCgr{p} (K,L)$.
The \emph{$p$-th (reduced) relative cohomology group} is
$\corHgr{p} (K,L) = \coZgr{p} (K,L) / \coBgr{p} (K, L)$.
For example, if $L = \emptyset$, then $\corHgr{p} (K,L) = \corHgr{p} (K)$,
and if $L = \{\emptyset\}$, then $\corHgr{p} (K,L)$ is isomorphic
to the conventional cohomology group in which $K$ does not contain the empty cell.
For relative cohomology, we write
$\corBetti{p}(K,L) = \rank{\corHgr{p}(K,L)} = \log_2 \card{\corHgr{p}(K,L)}$.

\ourparagraph{Matrix reduction.}
The classic algorithm for homology and cohomology reduces the boundary matrix
to Smith normal form; see \cite[$\S 11$]{Mun84}.
For modulo-$2$ arithmetic, this simplifies to Gaussian elimination.
We introduce special versions of this algorithm that are easy
to relate to the pertinent algebraic information, including the ranks of the groups.
To compute homology, we initialize $R = \partial$ and $U = \mbox{\rm Id}$
and reduce $R$ using left-to-right column additions
while maintaining the relation $R = \partial U$.
Write $\Low{j}$ for the row index of the lowest non-zero item in column $j$ of $R$,
and set $\Low{j} = -\infty$ if the column is zero.
\begin{tabbing}
mn\=mn\=mn\=mn\=mn\=mn\=mn\=mn\=\kill
{\tt Exhaustive column reduction algorithm:}                              \\*
\> {\tt for} $j=0$ {\tt to} $m$ {\tt do}                                  \\*
\> \> {\tt while} $\exists \ell<j$ with $\Low{\ell} > -\infty$
                                    and $R[\Low{\ell},j] \neq 0$ {\tt do} \\*
\> \> \> $R[.,j] = R[.,j] + R[.,\ell]$;
         $U[.,j] = U[.,j] + U[.,\ell]$.
\end{tabbing}
We call this algorithm \emph{exhaustive} because it attempts to remove non-zero
entries in column $j$ even after the lowest such entry has been established.
While this strategy has also been used in \cite{EdZo03},
it is different from the standard reduction algorithm used in persistent homology,
which proceeds to column $j+1$ as soon as $\Low{j}$ is established.
An important difference is that for the exhaustive reduction algorithm,
the produced matrices $R$ and $U$ can be uniquely defined in terms of
their algebraic structure,
so they do not depend on the algorithm that computes them.
In contrast, the standard reduction algorithm returns matrices $R$ and $U$
that are generally not unique and depend on the chosen order of column additions.
To relate the algorithm to the homology group of $K$, we interpret
the reduction of column $j$ as adding $\ssx_j$ to the complex.
Since we assume a monotonic ordering, we have a polyhedral complex after
every addition.
There are two possible outcomes when we add $\ssx_j$ with $\dime{\ssx_j} = p$:
\medskip \begin{itemize}
  \item  column $j$ is reduced to zero, in which case $\rBetti{p}$ increases by $1$;
  \item  column $j$ remains non-zero, in which case $\rBetti{p-1}$ decreases by $1$;
\end{itemize} \medskip
see \cite{DeEd95}.
In the first case, we say $\ssx_j$ \emph{gives birth} to a $p$-cycle,
and in the second case, we say $\ssx_j$ \emph{gives death} to a $(p-1)$-cycle,
namely the one given birth to by $\ssx_i$ with $i = \Low{j}$;
see \cite[Chapter VII]{EdHa10}.
At completion, $\rBetti{p}$ is the number of $p$-cells, $\ssx_j$,
such that column $j$ of $R$ is zero and $j \neq \Low{\ell}$ for all $\ell$.
Writing $\nSxB{p}$ for the number of $p$-cells that give birth
and $\nSxD{p}$ for the number that give death,
we have $\nSx{p} = \nSxB{p} + \nSxD{p}$ and $\rBetti{p} = \nSxB{p} - \nSxD{p+1}$.
We can therefore express $\nSxB{p}$ and $\nSxD{p}$ in terms of
the $\nSx{q}$ and the $\rBetti{q}$,
and since the Betti numbers are topological invariants,
we conclude that the $\nSxB{p}$ and the $\nSxD{p}$ neither depend on the
particular reduction algorithm nor on the ordering of the cells.

To compute cohomology, we initialize $Q = \partial$ and $V = \mbox{\rm Id}$
and reduce $Q$ using bottom-to-top row operations while maintaining the relation
$Q = V \partial$.
Write $\Left{i}$ for the column index of the leftmost non-zero entry
in row $i$ of $Q$, and set $\Left{i} = \infty$ if the row is zero.
\begin{tabbing}
mn\=mn\=mn\=mn\=mn\=mn\=mn\=mn\=\kill
{\tt Exhaustive row reduction algorithm:}                                  \\*
\> {\tt for} $i=m$ {\tt downto} $0$ {\tt do}                               \\*
\> \> {\tt while} $\exists \ell>i$ with $\Left{\ell} < \infty$
                                    and $Q[i,\Left{\ell}] \neq 0$ {\tt do} \\*
\> \> \> $Q[i,.] = Q[i,.] + Q[\ell,.]$;
         $V[i,.] = V[i,.] + V[\ell,.]$.
\end{tabbing}
Similar to before, we distinguish between two possible outcomes:
that a row is reduced to zero and that it remains non-zero.
To relate the algorithm to cohomology, 
we interpret the state of $Q$ after reducing $\ssx_i$ as the relative cohomology
of the pair, $(K, L)$, in which $L \subseteq K$ consists of
cells $\ssx_0$ to $\ssx_{i-1}$.
The reduction of a row is therefore akin to moving a cell
from $L$ to $K \setminus L$.
The two possible outcomes correspond again to \emph{births} and \emph{deaths},
this time of relative cocycles.
Writing $\ncoSxB{p}$ and $\ncoSxD{p}$ for the numbers of $p$-cells
of the two types, we have
$\ncoSx{p} = \ncoSxB{p} + \ncoSxD{p}$ and $\corBetti{p} = \ncoSxB{p} - \ncoSxD{p-1}$,
and as before these numbers neither depend on the
particular reduction algorithm nor on the ordering of the cells.

\ourparagraph{Duality.}
The births and deaths recorded in $R$ and in $Q$ are not the same
but they are closely related.
This is not surprising since a classic result in algebraic topology
asserts that the $\Zspace / 2 \Zspace$ ranks of the homology
and the cohomology groups coincide;
see e.g.\ \cite{Hat02,Mun84}.
We formally state this result for later reference,
and we relate it to the above reduction algorithms.
\begin{proposition}[Duality]
  \label{prop:Duality}
  Every polyhedral complex satisfies $\rBetti{p} = \corBetti{p}$,
  for all dimensions $p$.
\end{proposition}
\ourproof
  We compare the ranks, which we read off the reduced boundary matrices, $R$ and $Q$.
  As implied by a more general theorem in \cite{CEM06},
  $R$ and $Q$ have the same birth-death pairs; that is:
  $i = \Low{j}$ in $R$ iff $j = \Left{i}$ in $Q$.
  The number of birth-death pairs $(\ssx_i, \ssx_j)$ with
  $\dime{\ssx_i} = \dime{\ssx_j} - 1 = p-1$ is $\nSxD{p} = \ncoSxD{p-1}$.
  The $p$-th reduced Betti number is
  $\rBetti{p} = \nSxB{p} - \nSxD{p+1} = \nSx{p} - \nSxD{p} - \nSxD{p+1}$,
  and the rank of the $p$-th reduced cohomology group is
  $\corBetti{p} = \ncoSxB{p} - \ncoSxD{p-1} = \ncoSx{p} - \ncoSxD{p} - \ncoSxD{p-1}$,
  which we can now see are equal.
\eop

\section{Tri-partition}
\label{sec:3}

This section presents the first result of this paper:  the tri-partition
of a polyhedral complex in which the three sets represent unique aspects
of the complex' topology.
We begin with the introduction of the sets.

\ourparagraph{Trees, cotrees, and else.}
Letting $K$ be a polyhedral complex, we recall that a $p$-chain is a subset
of its $p$-cells.
A \emph{$p$-tree} is a $p$-chain, $A_p \subseteq K^p$,
that does not contain any non-empty $p$-cycle;
compare with the definition of a generalized tree in \cite{Kal83}.
Sometimes these generalized trees are referred to as \emph{acycles} \cite{STY17},
which motivates our notation.
A $p$-tree is \emph{maximal} if it is not properly contained in another $p$-tree.
Similarly, a \emph{$p$-cotree} is a $p$-cochain, $A^p \subseteq K^p$,
that does not contain any non-empty $p$-cocycle,
and it is \emph{maximal} if it is not properly contained in another $p$-cotree.
As examples consider the complexes in Figure \ref{fig:disk}.
Adding the triangles, quadrangles, and the outer face to
the graph in the left panel, we get a $2$-dimensional complex with $32$ edges.
Half the edges form a maximal $1$-tree, with the other half forming
a maximal $1$-cotree.
Moving from the left to the right panel, we get the annulus
by removing the center vertex
together with the incident edges and triangles as well as the outer region.
The remaining $24$ edges contain a maximal $1$-tree of size $15$ and
a maximal $1$-cotree of size $8$, leaving one edge unused.

Our sole requirement for the third set of $p$-cells, $E_p$,
is that its cardinality be $\card{E_p} = \rBetti{p}$.
Since we talk about partitions, we have $A_p \cap A^p = \emptyset$
and $E_p = K^p \setminus A_p \setminus A^p$,
which we will see implies the existence of $\rBetti{p}$ $p$-cycles that generate
$\rHgr{p}$ and of $\corBetti{p}$ $p$-cocycles that generate $\corHgr{p}$
such that each $p$-cell in $E_p$ belongs to exactly one of these
cycles and to exactly one of these cocycles;
see Section \ref{sec:4} for details.

\ourparagraph{Statement and proof.}
We give a constructive proof that $K^p$ permits a tri-partition as described.
More specifically, we construct such a tri-partition for every ordering
of the $p$-cells.
The ordering of the other cells is not important as long as the
overall ordering of $K$ is monotonic.
\begin{theorem}[Tri-partition]
  \label{thm:TriPartition}
  Let $K$ be a polyhedral complex.
  Then there exist tri-partitions $A_p \sqcup A^p \sqcup E_p = K^p$,
  for every dimension $p$,
  such that $A_p$ is a maximal $p$-tree, $A^p$ is a maximal $p$-cotree,
  and $E_p = K^p \setminus A_p \setminus A^p$ with $\card{E_p} = \rBetti{p}$.
\end{theorem}

It is possible to argue the existence of the
tri-partition in terms of column- and row-spaces of the boundary matrix.
Our proof is along the same lines but more specific and designed
to reveal additional properties that will be exploited
in the construction of bases in Section \ref{sec:4}.
\ourproof
  We get the tri-partition in three steps:
  assuming a fixed monotonic ordering of $K$,
  we first construct $A_p$, we second construct $A^p$,
  and we let $E_p$ contain the remaining $p$-cells.

  To get started, we sort the rows and columns of the boundary matrix
  according to the monotonic ordering of $K$.
  For the first step, we use the exhaustive column reduction algorithm formally stated
  in Section \ref{sec:2}.
  Proceeding from left to right,
  column $j$ is either a combination of preceding columns,
  in which case it gets reduced to zero,
  or it is independent of the preceding columns,
  in which case it remains non-zero.
  In the latter case, we add the $j$-th cell to $A_p$,
  in which $p$ is the dimension of this cell.
  At termination, $A_p$ is a maximal $p$-tree and $\card{A_p} = \nSxD{p}$
  by construction.
  For the second step, we use the exhaustive row reduction algorithm also formally stated
  in Section \ref{sec:2}.
  Proceeding from bottom to top, row $i$ is either a combination of
  succeeding (lower) rows, in which case it gets reduced to zero,
  or it is independent of the succeeding rows, in which case it remains non-zero.
  In the latter case, we add the $i$-th cell to $A^p$,
  in which $p$ is the dimension of this cell.
  At termination, $A^p$ is a maximal $p$-cotree and $\card{A^p} = \ncoSxD{p}$
  by construction.

  It remains to prove that $A_p \cap A^p = \emptyset$ and
  that $E_p = K^p \setminus A_p \setminus A^p$ has cardinality $\rBetti{p}$.
  To prove disjointness, we note that $k = \Low{\ell}$ after column reduction
  iff $\ell = \Left{k}$ after row reduction; see \cite{CEM06}.
  Writing $\ssx_i$ for the $i$-th cell in the monotonic ordering,
  we have $\ssx_i \in A^p$ iff $p = \dime{\ssx_i}$ and $\Left{i} < \infty$.
  Writing $j = \Left{i}$, this is equivalent to $i = \Low{j}$,
  which implies that $\ssx_i$ gives birth to the $p$-cycle that $\ssx_j$ destroys.
  Hence, $\ssx_i \not\in A_p$, as desired.
  The symmetric argument shows that $A^p$ contains no cells of $A_p$,
  which implies $A_p \cap A^p = \emptyset$.
  Setting $E_p = K^p \setminus A_p \setminus A^p$,
  we observe that it contains
  a $p$-cell iff neither the corresponding row nor the corresponding column
  contains a birth-death pair.
  In other words, each such $p$-cell gives birth to an essential $p$-cycle
  in homology and, equivalently,
  it gives birth to an essential $p$-cocycle in cohomology.
  There are $\rBetti{p} = \corBetti{p}$ of each kind,
  hence $\card{E_p} = \rBetti{p}$, as claimed.
\eop

The proof shows slightly more than claimed in Theorem \ref{thm:TriPartition},
namely that there is a unique tri-partition for every monotonic ordering of $K$.
On the other hand, two different monotonic orderings do not necessarily
have different tri-partitions.
For example, the tri-partition in dimension $p$ is invariant as long as
we retain the ordering among the $p$-cells,
rearranging the other cells at will provided the overall ordering
remain monotonic.
This suggests we consider the collection of tri-partitions generated
by monotonic orderings of $K$.
Looking at its three constituents,
we note that the collections of sets $A_p$,
of sets $A^p$, and of sets $E_p$ are three matroids;
see the formal claim and the proof in Appendix \ref{app:A}.

\begin{figure}[hbt]
  \centering \resizebox{!}{2.2in}{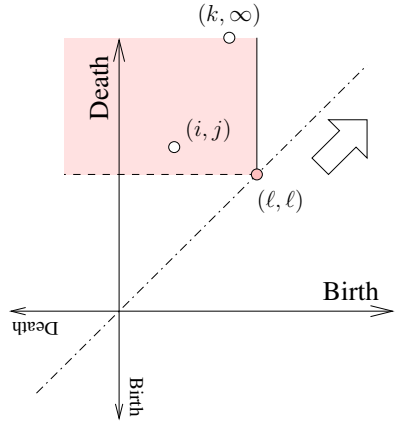}
  \caption{The persistence diagram of a monotonic ordering reveals
    the ranks of the reduced homology groups of all complexes $K_\ell$,
    and after reflecting the points at infinite from top to left,
    it reveals also the ranks of the reduced relative cohomology groups of all pairs
    $(K, K_\ell)$.}
  \label{fig:diagram}
\end{figure}

\ourparagraph{Tri-partitions and persistence diagrams.}
It is interesting to compare the tri-partition with the persistence
diagram for the same monotonic ordering.
Let $\ssx_0, \ssx_1, \ldots, \ssx_m$ be such an ordering
and write $K_\ell = \{ \ssx_0, \ssx_1, \ldots, \ssx_\ell \}$
for every $0 \leq \ell \leq m$.
The persistence diagram consists of all points $(i,j)$ for which
$R[.,i] = 0$ and $i = \Low{j}$
and all points $(k, \infty)$ for which $R[.,k] = 0$
but $k \neq \Low{\ell}$ for all $0 \leq \ell \leq m$;
see Figure \ref{fig:diagram} and refer to \cite{EdHa10} for details.
Importantly, the number of points in upper-left quadrants anchored
at points on the diagonal give the reduced Betti numbers of
complexes in the filtration.
Specifically, $\rBetti{p} (K_\ell)$
is the number of points $(i,j)$ with $\dime{\ssx_i} = p$
that satisfy $i \leq \ell < j$,
which includes the case $j = \infty$.
As we slide the quadrant to the right and up the diagonal,
we can read the reduced Betti numbers of all complexes in the filtration.
With a few modifications, we can also read the ranks
of the relative cohomology groups:
reverse the two axes and exchange birth with death,
move the points at infinity from north to west by reflecting
them across the minor diagonal,
and exchange the closed and open sides of the quadrant,
which we slide to the left and down the diagonal.

The persistence diagram implies the tri-partition
but not the other way round.
Specifically, for every finite point $(i,j)$
with $\dime{\ssx_i} = p$ in the diagram,
we have $\ssx_i \in A^p$ and $\ssx_j \in A_{p+1}$,
and for every point $(k, \infty)$ with $\dime{\ssx_k} = p$
we have $\ssx_k \in E_p$.
In other words,
the tri-partition records which cells give birth, which of those
are essential, and which cells give death,
but it does not determine the pairing that defines the persistence diagram.

\ourparagraph{Incremental construction.}
We conclude this section with a brief discussion of the incremental construction
of the tri-partition.
Suppose we have the tri-partition of $K_{\ell}$,
how can we modify it to get the tri-partition of $K_{\ell + 1}$?
There are only two cases, depending on whether $R[.,\ell+1] = 0$
after column reduction or not.
Let $p = \dime{\ssx_{\ell+1}}$.
\medskip \begin{description}
  \item[{\sc Case}] $R[.,\ell+1]    = 0.$
    Then $\ssx_{\ell+1}$ gives birth to a $p$-cycle,
    so we add $\ssx_{\ell+1}$ to $E_p$, leaving $A_p$ and $A^p$ untouched.
  \item[{\sc Case}] $R[.,\ell+1] \neq 0.$
    Then $\ssx_{\ell+1}$ gives death to a $(p-1)$-cycle,
    and we add $\ssx_{\ell+1}$ to $A_p$.
    Letting $k = \Low{\ell+1}$, we have $\ssx_k \in E_{p-1}$
    and since its class just got killed, we move it to $A^{p-1}$.
\end{description} \medskip
Note the asymmetry between the trees and the cotrees revealed by the incremental construction.
Particularly perplexing, at first, is the move of $\ssx_k \in E_{p-1}$
--- which signifies a birth in relative cohomology ---
to $\ssx_k \in A^{p-1}$ --- which signifies a death in relative cohomology.
The reason for this drastic change is of course the difference in direction,
which is from left to right in the incremental construction,
and from right to left in the computation of relative cohomology.

\section{Bases}
\label{sec:4}

Besides constructing tri-partitions of a polyhedral complex,
the exhaustive reduction algorithms compute canonical bases in homology
and in cohomology.
This section describes these bases and proves some of their properties.

\ourparagraph{Cycles and chains.}
Fixing a monotonic ordering of a polyhedral complex, $K$,
we write $K^p = A_p \sqcup A^p \sqcup E_p$ for the corresponding tri-partition
of the $p$-cells.
For each $\ssx_j \in K^p$, we define a unique $p$-cycle or a unique $p$-chain
with non-empty boundary.
Specifically, if $\ssx_j \in A^p \sqcup E_p$, then there is a unique $p$-cycle
$\Cycle{p}{\ssx_j} \subseteq A_p \sqcup \{\ssx_j\}$,
which we refer to as the \emph{canonical $p$-cycle} of $\ssx_j$.
If $\ssx_j \in A_p$, then there is a unique $p$-chain,
$\Chain{p}{\ssx_j} \subseteq A_p$, whose boundary is the
sum of canonical $(p-1)$-cycles that is rendered trivial by the
addition of $\ssx_j$ to $K_{j-1}$,
which we refer to as the \emph{canonical $p$-chain} of $\ssx_j$.
Symmetrically, for every $\ssx_i \in K^p$, we define the
\emph{canonical $p$-cocycle},
$\Cocycle{p}{\ssx_i} \subseteq A^p \sqcup \{\ssx_i\}$,
if $\ssx_i \in A_p \sqcup E_p$,
and the \emph{canonical $p$-cochain}, $\Cochain{p}{\ssx_i} \subseteq A^p$,
if $\ssx_i \in A^p$.
We prove a technical lemma.
\begin{lemma}[Off-diagonal Entries]
  \label{lem:OffdiagonalEntries}
  Let $R = \partial U$ and $Q = V \partial$ be the matrix equations after
  exhaustive reduction.
  For every $i \neq j$ there is a dimension $q$ such that
  $U[i,j] = 1$ implies $\ssx_i \in A_q$
  and $V[i,j] = 1$ implies $\ssx_j \in A^q$.
\end{lemma}
\ourproof
  To prove that all non-zero off-diagonal entries in $U$ belong to rows
  of cells in $A_q$, for some $q$, we note that this
  is trivially true at the start of the reduction algorithm,
  when $U = \mbox{\rm Id}$.
  A column $\ell$ is added to column $j$ only if $\ell < j$
  and $R[.,\ell]$ is non-zero.
  Since the algorithm proceeds from left to right,
  this implies that $\ssx_\ell \in A_q$.
  Assuming inductively that also all off-diagonal non-zero entries
  in $U[.,\ell]$ belong to rows
  of $q$-cells in $A_q$, we see that the column operation maintains the
  claim about off-diagonal entries.

  The argument why all non-zero off-diagonal entries in $V$ belong to
  columns of cells in $A^q$ is symmetric and omitted.
\eop

The technical lemma is useful to shed light on the connection between
the canonical cycles, chains, cocycles, cochains and the matrices
after exhaustive reduction.
\begin{lemma}[Columns and Rows]
  \label{lem:ColumnsNRows}
  After exhaustive column reduction of $R = \partial U$,
  the columns of $U$ store the canonical cycles and chains,
  and after exhaustive row reduction of $Q = V \partial$,
  the rows of $V$ store the canonical cocycles and cochains:
  \begin{align}
    U[.,j]  &=  \left\{ \begin{array}{ll}
                  \Cycle{p}{\ssx_j}  &  \mbox{\rm if~} \ssx_j \in A^p \sqcup E_p , \\
                  \Chain{p}{\ssx_j}  &  \mbox{\rm if~} \ssx_j \in A_p ,
                \end{array} \right.           \\
    V[i,.]  &=  \left\{ \begin{array}{ll}
                  \Cocycle{p}{\ssx_i}  &  \mbox{\rm if~} \ssx_i \in A_p \sqcup E_p , \\
                  \Cochain{p}{\ssx_i}  &  \mbox{\rm if~} \ssx_i \in A^p .
                \end{array} \right.          
  \end{align}
\end{lemma}
\ourproof
  Because of symmetry, it suffices to prove the claims about the cycles and chains.
  Consider first the case in which $\ssx_j \in A^p \sqcup E_p$.
  After completing the reduction of column~$j$, $U[.,j]$ stores a cycle.
  All cells in this cycle have the same dimension as $\ssx_j$, which is $p$.
  Lemma~\ref{lem:OffdiagonalEntries} implies that
  this cycle is a subset of $A_p \sqcup \{\ssx_j\}$.
  There is only one such cycle, namely $\Cycle{p}{\ssx_j}$,
  which implies that $U[.,j]$ stores this cycle, as claimed.

  Consider second the case in which $\ssx_j \in A_p$.
  To show that $U[.,j]$ stores $\Chain{p}{\ssx_j}$, we note that all non-zero
  entries in column $j$ of $U$ belong to rows of cells in $A_p$,
  and this includes the diagonal entry.
  Writing $\cgr{} \subseteq A_p$ for this chain and $\zgr{} = \partial \cgr{}$
  for its boundary, we note that $\zgr{}$ is stored in column $j$ of $R$.
  Separating the birth-giving from the death-giving $(p-1)$-cells,
  we write $\zgr{} = \zgr{\rm bth} \sqcup \zgr{\rm dth}$.
  Note that $\zgr{\rm bth} \subseteq \zgr{} \subseteq A_{p-1} \sqcup \zgr{\rm bth}$
  and that there is only one such $(p-1)$-cycle, namely the sum of the
  canonical $(p-1)$-cycles of the $\ssx \in \zgr{\rm bth}$.
  Hence, $\zgr{} = \sum_{\ssx \in \zgr{\rm bth}} \Cycle{p-1}{\ssx}$.
  Since the column reduction algorithm is exhaustive, each $(p-1)$-cycle
  in the sum is born before $\ssx_j$ and dies after $\ssx_j$.
  Any other sum of non-trivial canonical $(p-1)$-cycles is non-homologous to $\zgr{}$.
  By construction, $\zgr{}$ goes from non-trivial to trivial when we add $\ssx_j$,
  which implies that $U[.,j]$ stores $\Chain{p}{\ssx_j}$, as claimed.
\eop

\ourparagraph{Canonical bases.}
The columns of $U$ and $R$ provide bases for the cycle, the boundary, and the
homology groups, and the rows of $V$ and $Q$ provide bases for the cocycle,
the coboundary, and the cohomology groups.
These bases depend on the ordering of the cells, but they are \emph{canonical}
in the sense that they are defined in terms of their algebraic properties
and do not depend on the algorithms that compute them.
\begin{theorem}[Canonical Bases]
  \label{thm:CanonicalBases}
  Assume a monotonic ordering of a polyhedral complex, $K$,
  and let $K^p = A_p \sqcup A^p \sqcup E_p$ be the corresponding tri-partition.
  Then
  \begin{itemize}\denselist
    \item $\{ \Cycle{p}{\ssx_j} \mid \ssx_j \in A^p \sqcup E_p \}$
      is a basis of $\Zgr{p} (K)$.
    \item $\{ \Cycle{p}{\ssx_j} \mid \ssx_j \in E_p \}$
      generates a basis of $\rHgr{p} (K)$.
    \item $\{ \partial \Chain{p}{\ssx_j} \mid \ssx_j \in A_p \}$
      is a basis of $\Bgr{p-1} (K)$.
    \item $\{ \Cocycle{p}{\ssx_i} \mid \ssx_i \in A_p \sqcup E_p \}$
      is a basis of $\coZgr{p} (K)$.
    \item $\{ \Cocycle{p}{\ssx_i} \mid \ssx_i \in E_p \}$
      generates a basis of $\corHgr{p} (K)$.
    \item $\{ \delta \Cochain{p}{\ssx_i} \mid \ssx_i \in A^p \}$
      is a basis of $\coBgr{p+1} (K)$.
  \end{itemize}
\end{theorem}
\ourproof
  Because of symmetry, we can limit ourselves to the first three claims,
  which are about cycles and chains.
  We prove these claims in sequence.

  To see that the $\Cycle{p}{\ssx_j}$, over all $\ssx_j \in A^p \sqcup E_p$,
  form a basis of $\Zgr{p} (K)$,
  we note that these cycles are clearly independent.
  Let $\zgr{}$ be an arbitrary $p$-cycle, and write $\zgr{j} = \Cycle{p}{\ssx_j}$
  for every $\ssx_j \in (A^p \sqcup E_p) \cap \zgr{}$.
  Then $\zgr{} = \sum_j \zgr{j}$, for if they were different,
  then $\zgr{} + \sum_j \zgr{j}$ would be a $p$-cycle contained in $A_p$,
  which contradicts the acyclicity of $A_p$.

  Recall that a homology class is essential if it is non-trivial in $K$.
  By construction, when $\ssx_j \in E_p$, then $\Cycle{p}{\ssx_j}$ generates an
  essential class, and when $\ssx_j \in A^p$,
  then $\Cycle{p}{\ssx_j}$ is trivial or homologous to a sum of cycles
  defined by $p$-cells in $E_p$.
  Since $\rHgr{p} (K)$ requires $\rBetti{p}$ generators and there are only
  $\rBetti{p}$ cells in $E_p$, each homology class represented by such a cell
  must be a generator of $\rHgr{p} (K)$.

  To see that the $\partial \Chain{p}{\ssx_j}$, over all $\ssx_j \in A_p$,
  form a basis of $\Bgr{p-1} (K)$, we note that these cycles are independent.
  Indeed, if they were not independent, then we had a non-empty sum of chains
  with empty boundary, which contradicts that the chains are all contained
  in $A_p$ since $A_p$ contains no cycle by construction.
  To show that the $\partial \Chain{p}{\ssx_j}$ span the $(p-1)$-dimensional
  boundary group, we recall that the number of $(p-1)$-cycles
  $\partial \Chain{p}{\ssx_j}$ is $\card{A_p} = \nSxD{p}$.
  For comparison, the rank of ${\Bgr{p-1} (K)}$ is equal to the rank of $\Zgr{p-1} (K)$
  minus the rank of $\rHgr{p-1} (K)$,
  which is $\nSxB{p-1} - (\nSxB{p-1} - \nSxD{p}) = \nSxD{p}$.
  Since this is the same as the number of $(p-1)$-cycles,
  we conclude that the $\partial \Chain{p}{\ssx_j}$ indeed form a basis of $\Bgr{p-1}$.
\eop

Theorem \ref{thm:CanonicalBases} implies the algebraic analog of the
Helmholtz--Hodge decomposition, namely that the $p$-th chain group satisfies
$\Cgr{p} \simeq \Bgr{p-1} \oplus \coBgr{p+1} \oplus \rHgr{p}$
for every dimension $p$.
This is the algebraic way of saying that each $p$-cell
either kills a $(p-1)$-cycle, gives birth to a $p$-cycle that later dies,
or gives birth to an essential $p$-cycle.
Indeed, if we construct the filtration in reverse while maintaining
the relative cohomology, the second of these three options correspond
to killing a $(p+1)$-cocycle.
To get the algebraic decomposition in standard form,
$\Cgr{p} \simeq \Bgr{p} \oplus \coBgr{p} \oplus \rHgr{p}$,
we note $\Bgr{p-1} \simeq \coBgr{p}$ and $\coBgr{p+1} \simeq \Bgr{p}$ as needed.

\ourparagraph{Intersections of basis vectors.}
To study the relation between the various basis vectors,
we consider the matrix product, $VU$,
which we compute over $\Zspace$ so that $1+1=2$.
To predict its entries, we begin with the special case in which
the column of $U$ stores a dead cycle and the row of $V$
stores a dead cocycle.
\begin{lemma}[Two Crossings]
  \label{lem:TwoCrossings}
  Assume a monotonic ordering of a polyhedral complex, $K$,
  let $K^p = A_p \sqcup A^p \sqcup E_p$ be the corresponding tri-partition,
  and suppose $\ssx_i \in A_p$ and $\ssx_j \in A^p$.
  Then $\ssx_i \in \Cycle{p}{\ssx_j}$ iff $\ssx_j \in \Cocycle{p}{\ssx_i}$.
\end{lemma}
\ourproof
  Recall that $\ssx_j \in A^p$ implies that $\ssx_j$ gives birth to a $p$-cycle
  that dies before the filtration ends.
  The $p$-cycle born when we add $\ssx_j$ to $K_{j-1}$ is $\Cycle{p}{\ssx_j}$,
  and we write $\cgr{p+1}$ for the $(p+1)$-chain that gives death to this $p$-cycle.
  The death occurs either because the cycle becomes trivial or it becomes homologous
  to another cycle.

  Consider first the easier case, when $\partial \cgr{p+1} = \Cycle{p}{\ssx_j}$.
  Suppose $\ssx_i \in \Cycle{p}{\ssx_j}$.
  Since $\Cocycle{p}{\ssx_i}$ is a cocycle, every $(p+1)$-cell in $\cgr{p+1}$
  has an even number of $p$-faces in $\Cocycle{p}{\ssx_i}$.
  The sum of these even numbers, over all $(p+1)$-cells in $\cgr{p+1}$,
  is of course even,
  and this even number is also the sum over all $p$-cells in $\Cocycle{p}{\ssx_i}$
  of their numbers of cofaces in $\cgr{p+1}$.
  This number of cofaces is odd for every $p$-cell in $\Cycle{p}{\ssx_j}$,
  and even for every other $p$-cell.
  Now $\ssx_j$ is the only $p$-cell in $\Cycle{p}{\ssx_j}$ that does not belong to $A_p$,
  which implies that $\ssx_j$ is the only $p$-cell of $\Cycle{p}{\ssx_j}$
  --- other than $\ssx_i$ --- that is possibly in $\Cocycle{p}{\ssx_i}$.
  Since $\ssx_i$ has an odd number of cofaces in $\cgr{p+1}$ and the sum is even,
  this implies that $\ssx_j \in \Cocycle{p}{\ssx_i}$,
  as desired.

  Consider second the case in which
  $\partial \cgr{p+1} = \Cycle{p}{\ssx_j} + \sum_{k \in \Lambda} \Cycle{p}{\ssx_{k}}$
  for a non-empty index set $\Lambda$.
  By construction, $k < j$ for every $k \in \Lambda$.
  This implies $\ssx_{k} \in E_p$ at time $j$.
  To get a contradiction, we assume $\ssx_{k} \in \Cocycle{p}{\ssx_i}$.
  Recall that $\ssx_{k} \in E_p$ implies that $\Cocycle{p}{\ssx_{k}}$
  has empty intersection with $A_p$, so $\ssx_i \not\in \Cocycle{p}{\ssx_{k}}$.
  Let $\cozgr{p} = \Cocycle{p}{\ssx_i} + \Cocycle{p}{\ssx_{k}}$,
  which is again a $p$-cocycle that contains $\ssx_i$.
  But now we have two $p$-cocycles in $A^p \sqcup \{\ssx_i\}$,
  which contradicts that $A^p$ is a cotree.
  Repeating the argument of the first case,
  we conclude that $\ssx_j \in \Cocycle{p}{\ssx_i}$
  since none of the $\ssx_k$ belongs to $\Cocycle{p}{\ssx_i}$, as desired.
\eop

\begin{figure}[hbt]
    \centering \resizebox{!}{1.6in}{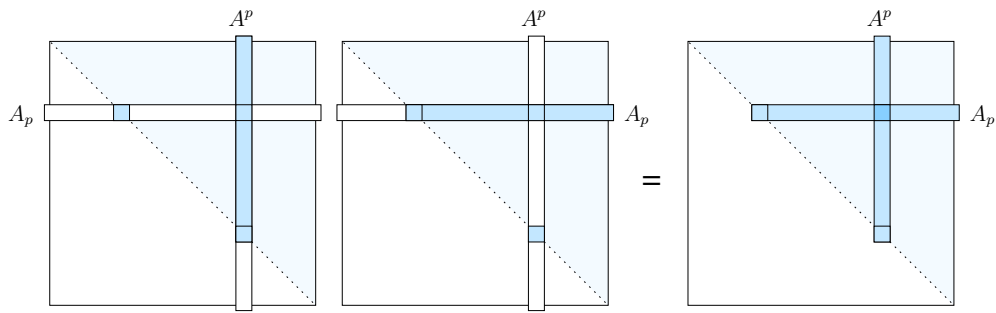}
    \caption{We get $VU$ by multiplying $V$
      --- whose rows contain the basis vectors in cohomology ---
      with $U$
      --- whose columns contain the basis vectors in homology.
      To compute the off-diagonal entries of $VU$,
      only the columns of cells in $A^p$ of $V$
      and the rows of cells in $A_p$ of $U$ are relevant.}
    \label{fig:VU}
  \end{figure}

The above lemma covers only one of the nine possible combinations.
In each case, the number of cells in which the chain and the cochain
overlap is either $0$, $1$, or $2$.
To formulate this claim in greater detail, we write $\xgr{p}$ for $\cgr{p}$ or $\zgr{p}$
and $\coxgr{p}$ for $\cocgr{p}$ and $\cozgr{p}$ in cases in which the parameter
decides which of the two functions applies.
\begin{theorem}[Intersection Patterns]
  \label{thm:IntersectionPatterns}
  Assume a monotonic ordering of a polyhedral complex, $K$,
  and let $K^p = A_p \sqcup A^p \sqcup E_p$ be the corresponding tri-partition
  in dimension $p$.
  Then
  \begin{align}
    VU [i,j]  &=  \left\{ \begin{array}{ll}
                    2  &  \mbox{\rm if~}  (\ssx_i, \ssx_j) \in A_p \times A^p
                                          \mbox{\rm ~and~} \ssx_i \in \Cycle{p}{\ssx_j};   \\
                    1  &  \mbox{\rm if~}  $i=j$
                          \mbox{\rm ~or~} (\ssx_i, \ssx_j) \in A_p \times (A_p \sqcup E_p)
                                          \mbox{\rm ~and~} \ssx_i \in \Xain{p}{\ssx_j}     \\
                       &  ~~~~~~~
                          \mbox{\rm ~or~} (\ssx_i, \ssx_j) \in (A^p \sqcup E_p) \times A^p
                                          \mbox{\rm ~and~} \ssx_j \in \Coxain{p}{\ssx_i};  \\
                    0  &  \mbox{\rm otherwise.}
                  \end{array} \right.
    \label{eqn:pattern}
  \end{align}
\end{theorem}
\ourproof
  We illustrate the argument in Figure \ref{fig:VU},
  which shows $V$, $U$, $VU$ together with the rows for one cell in $A_p$
  and the columns for one cell in $A^p$.
  Recall that $U$ and $V$ are both upper-triangular,
  with all diagonal entries equal to $1$.
  It follows that $VU$ is upper-triangular,
  with all diagonal entries equal to $1$ as well.
  To prove \eqref{eqn:pattern} for the off-diagonal entries of $VU$,
  we recall Lemma \ref{lem:OffdiagonalEntries} and note that it implies
  \begin{align}
    VU [i,j]  &=  V[i,.] \cdot U[.,j]  =  V[i,i] \, U[i,j] + V[i,j] \, U[j,j] .
  \end{align}
  Indeed, if $\ssx_i \in A_p$, then $V[i,i] = 1$ is the only non-zero entry in row $i$
  that belongs to a column of a cell in $A_p$.
  Multiplying this row with $U$ has the effect of copying row $i$ of $U$
  to $VU$; see the shaded rows in Figure \ref{fig:VU}.
  Symmetrically, if $\ssx_j \in A^p$, then $U[j,j] = 1$ is the only
  non-zero entry in column $j$ that belongs to a row of a cell in $A^p$.
  Multiplying $V$ with this column has the effect of copying column $j$
  of $V$ to $VU$; see the shaded columns in Figure \ref{fig:VU}.

  The shaded rows and columns are added in $VU$,
  so the entries in their intersections can be $0$, $1$, or $2$.
  By Lemma \ref{lem:TwoCrossings}, they can only be $0$ or $2$,
  which implies the first line in \eqref{eqn:pattern}.
  The remaining entries in the shaded rows and columns of $VU$
  are as in $V$ and in $U$,
  and all other off-diagonal entries in $VU$ are zero,
  which implies the second and the third line in \eqref{eqn:pattern}.
\eop

Consider for example the case $(\ssx_i, \ssx_j) \in E_p \times E_p$.
Then Theorem \ref{thm:IntersectionPatterns} implies
$VU[i,j] = 1$ if $i=j$ and $VU[i,j] = 0$ if $i \neq j$.
In words, each $p$-cell of $E_p$ belongs to exactly one generating $p$-cycle
of $\rHgr{p}$ and to exactly one generating $p$-cocycle of $\corHgr{p}$,
and there are no other intersections between the basis vectors
of $\rHgr{p}$ and the basis vectors of $\corHgr{p}$.

\section{Discussion}
\label{sec:5}

The main contributions of this paper are the construction of
a tri-partition of a polyhedral complex and the analysis of the corresponding
bases in homology and in cohomology.
For a given monotonic ordering, the tri-partition is unique
and so are the corresponding bases.
We mention a few questions suggested by the work reported in this paper:
\medskip \begin{itemize}\denselist
  \item Our constructions generalize to situations in which
    homology and cohomology are defined for field coefficients.
    Do they also generalize to non-field coefficients, for example the integers?
  \item Can the analogy between the tri-partition and the
    Helmholtz--Hodge decomposition of a smooth vector field be used to gain
    insights on either side?
    For example, does the tri-partition lead to a fast algorithm for constructing
    harmonic cycles, that is, whose Laplacian is zero?
  \item Can the tri-partitions be used to shed light on the stochastic properties
    of simplicial complexes as studied in \cite{LiPe17}?
\end{itemize} \medskip
Applications of tri-partitions outside of mathematics are 
at least as important as finding connections within mathematics;
see \cite{EdOl18} for a first step.
Particularly interesting is the use of the trees and cotrees to explore
cave systems, such as within biomolecules and the molecular structure
of materials.

\subsubsection*{Acknowledgements}

\footnotesize{
The authors of this paper thanks J\"{o}rg Peters, Konrad Polthier,
and G\"{u}nter Rote for insightful discussions on the topic of this paper.
}

\newpage


\appendix \clearpage
\normalsize

\section{Matroids}
\label{app:A}

The fact that the greedy algorithm succeeds in constructing maximal trees
and maximal cotrees is not surprising since both form matroids \cite{Oxl92}.
We recall what this means.
Let $E$ be finite and $\Fcal$ a collection of subsets of $E$.
We call $(\Fcal, E)$ an \emph{abstract simplicial complex} if
$\emptyset \in \Fcal$ and $\Fcal$ is closed under taking subsets.
It is a \emph{matroid} if, in addition, $(\Fcal, E)$ satisfies
the \emph{exchange property}:  $F, G \in \Fcal$ with $\card{G} < \card{F}$
implies the existence of $e \in F$ such that $G \cup \{e\} \in \Fcal$.
Traditionally, the sets in $\Fcal$ are called \emph{independent},
and the exchange property implies that all maximal independent sets in $\Fcal$
have the same cardinality.
It is often convenient to focus on the maximal sets as all others are implied by inclusion.
The exchange property can be replaced by the following, equivalent property:
if $F, G$ are different maximal independent sets of $\Fcal$ and $a \in F \setminus G$, 
then there exists $b \in G \setminus F$ such that
$F \setminus \{a\} \cup \{b\} \in \Fcal$.

Given a polyhedral complex, $K$, we write $\Acal_p$, $\Acal^p$, and $\Ecal_p$
for the collections of $p$-trees, $p$-cotrees, and $p$-dimensional leftover,
which we recall are sets of $p$-cells.
We note that all these sets arise in tri-partitions of $K$ constructed for some
monotonic ordering.
Indeed, every $A_p \in \Acal_p$ arises as a subset of the maximal $p$-tree
if we order the $p$-cells in $A_p$ before all others.
Symmetrically, every $A^p \in \Acal^p$ arises as a subset of the maximal $p$-cotree
if we order the $p$-cells in $A^p$ after all others.
To see that every $E_p \in \Ecal_p$ arises as a leftover constructed
for some ordering,
we recall that there are collections
$D_p = \{ \Cycle{p}{\ssx} \mid \ssx \in E_p \}$ and
$D^p = \{ \Cocycle{p}{\ssx} \mid \ssx \in E_p \}$ that
satisfy Theorem \ref{thm:IntersectionPatterns};
see in particular the remark following the proof of this theorem.
Writing $\bigcup D_p$ and $\bigcup D^p$ for the $p$-cells that
belong to the cycles and cocycles in the two collections,
we note that $\bigcup D_p \setminus E_p$ and $\bigcup D^p \setminus E_p$
are disjoint, so we can make sure that the $p$-cells in
$\bigcup D_p \setminus E_p$ precede the $p$-cells in $E_p$,
and the latter precede the $p$-cells in $\bigcup D^p \setminus E_p$.
Adding the remaining $p$-cells arbitrarily, we get
$E_p$ as the leftover for the ordering.

It is well known that the trees and the cotrees have matroid structure.
We add that the same is true for the leftover sets.
\begin{lemma}[Tri-matroids]
  \label{lem:TriMatroids}
  Let $K$ be a polyhedral complex.
  Then $(\Acal_p, K^p)$, $(\Acal^p, K^p)$, and $(\Ecal_p, K^p)$
  are matroids for every dimension $p$.
\end{lemma}
\ourproof
  We prove the claim for $p$-trees as a warm-up exercise,
  skipping the argument for $p$-cotrees, which is almost verbatim the same.
  Let $F, G \in \Acal_p$ be maximal, and let $a \in F \setminus G$.
  Adding $a$ to $G$ creates a unique $p$-cycle, $A \subseteq G \cup \{a\}$.
  Adding a $p$-cell $b \in A \setminus F$ to $F$ creates again a unique $p$-cycle,
  which for the purpose of this proof we refer to as an \emph{elementary $p$-cycle}
  in $F \cup A$.
  The elementary $p$-cycles span the entire space of $p$-cycles of $F \cup A$,
  which includes $A$.
  We have $a \in A$, so $a$ must belong to at least one elementary $p$-cycle.
  Letting $b \in A \setminus F$ be a $p$-cell whose elementary $p$-cycle
  contains $a$, we get
  $F \setminus \{a\} \cup \{b\}$ as an independent set.
  Noting that $b \in G$, this implies that $( \Acal_p, K^p)$ is a matroid, as claimed.

  To prove that $(\Ecal_p, K^p)$ is a matroid,
  we use the fact that $(\Acal_p, K^p)$ and $(\Acal^p, K^p)$ are matroids,
  and that for each maximal $E_p \in \Ecal_p$ there are maximal $A_p \in \Acal_p$
  and $A^p \in \Acal^p$ such that $A_p \sqcup A^p \sqcup E_p = K^p$.
  Let $E_p' \in \Ecal_p$ be maximal and different from $E_p$,
  and let $A_p'$, ${A^p}'$ be a maximal $p$-tree and a maximal $p$-cotree
  with $A_p' \sqcup {A^p}' \sqcup E_p' = K^p$.
  Let $a \in E_p \setminus E_p'$ and assume without loss of generality
  that $a \in A_p'$.
  We add $a$ to $A_p$ and let $b \neq a$ be any $p$-cell of the thus created
  unique $p$-cycle.
  Removing $b$ from $A_p \cup \{a\}$, we get again a $p$-tree.
  If $b \in E_p'$, then we proceed to the next step,
  else $b \in {A^p}'$, we add $b$ to $A^p$,
  and we iterate with a $p$-cell $c$ in the thus created unique $p$-cocycle.
  Continuing this way, we eventually get a $p$-cell $z$ in $E_p'$.
  Indeed, every step makes $A_p$ more similar to $A_p'$
  or it makes $A^p$ more similar to ${A^p}'$, so the process must terminate.
  By construction, $z \not\in E_p$ and
  $E_p \setminus \{a\} \cup \{z\}$ is a maximal independent set of $\Ecal_p$,
  which implies that $(\Ecal_p, K^p)$ is a matroid, as claimed.
\eop

We note that the matroid structure implies that the greedy algorithm
can be used to construct optimal trees, cotrees, and leftovers efficiently.
This is in contrast to optimal bases, which for many objective functions
are NP-hard to construct \cite{ChFr11}.

The proof that $(\Ecal_p, K^p)$ is a matroid extends to general partitions
of a ground-set.
Fixing an integer $k$ and a set $E$,
we consider partitions $E = F_1 \sqcup F_2 \sqcup \ldots \sqcup F_k$
such that for each $1 \leq i < k$ the collection of sets $F_i$ is a matroid over $E$,
and conclude that the collection of sets $F_k$ is also a matroid over $E$.

\end{document}